\documentclass[11pt,a4paper,english]{article}

\usepackage{babel,latexsym,graphicx,amsmath,amsfonts,amsthm,amssymb}
\usepackage[colorlinks,citecolor=blue,urlcolor=blue,linkcolor=blue]{hyperref}
\usepackage[round]{natbib}
\usepackage{authblk}

\newtheoremstyle{localrem}
	{5pt} 
	{5pt} 
	{\rm} 
	{} 
	{\bf} 
	{{\rm.}} 
	{.7em} 
	{} 

\theoremstyle{localrem}
\newtheorem{Definition}{Definition}[section]

\newtheoremstyle{localthm}
	{5pt} 
	{5pt} 
	{\sl} 
	{} 
	{\bf} 
	{{\rm.}} 
	{.7em} 
	{} 

\theoremstyle{localthm}
\newtheorem{Corollary}[Definition]{Corollary}

\newtheorem{Lemma}[Definition]{Lemma}

\newcommand{\R}{\mathbb{R}}

\newcommand{\Rmda}{\R^m_{\downarrow}}
\newcommand{\Rmdx}{\R^m_{\downarrow, \bs{x}}}

\newcommand{\LL}{\mathcal{L}}
\newcommand{\MM}{\mathcal{M}}
\newcommand{\PP}{\mathcal{P}}
\newcommand{\UU}{\mathcal{U}}
\newcommand{\XX}{\mathcal{X}}
\newcommand{\WW}{\mathcal{W}}

\DeclareMathOperator*{\Corr}{Corr}
\DeclareMathOperator*{\argmin}{arg\,min}

\def\Pr{\mathop{\rm I\!P}\nolimits}

\def\bs{\boldsymbol}
\def\hat{\widehat}

\author[1]{Alexander Henzi\thanks{alexander.henzi@stat.unibe.ch}}
\author[2]{Alexandre M\"osching\thanks{alexandre.moesching@uni-goettingen.de}}
\author[1]{Lutz D\"umbgen\thanks{duembgen@stat.unibe.ch}}
\affil[1]{University of Bern, Switzerland}	
\affil[2]{Georg-August-University of G\"ottingen, Germany}
\date{\today}
\title{Accelerating the Pool-Adjacent-Violators Algorithm\\
	for Isotonic Distributional Regression}

\begin{document}

\maketitle

\begin{abstract}
In the context of estimating stochastically ordered distribution functions, the pool-adjacent-violators algorithm (PAVA) can be modified such that the computation times are reduced substantially. This is achieved by studying the dependence of antitonic weighted least squares fits on the response vector to be approximated.

\paragraph{Keywords:}
Monotone regression, sequential computation, weighted least squares

\paragraph{AMS 2000 subject classifications:}
62G08, 62G30, 62-08
\end{abstract}

\section{Introduction}
\label{sec:introduction}

Let $\mathcal{X}$ be a set equipped with a binary relation $\preceq$, for instance, some partial order. The general problem is as follows: For $m \ge 2$ pairs $(x_1,z_1),\ldots, (x_m,z_m) \in \XX\times \R$ and weights $w_1,\ldots,w_m > 0$, let 
\begin{equation}
\label{eq:A(z)}
	A(\bs{z}) \ := \ \argmin_{\bs{f} \in \Rmdx} \sum_{j=1}^m w_j (z_j - f_j)^2 ,
\end{equation}
where
\[
	\Rmdx
	\ := \
	\{\bs f\in\R^m: x_i\preceq x_j \ \text{implies that} \ f_i\ge f_j\}.
\]
Suppose that $\bs{z}^{(0)}, \bs{z}^{(1)}, \ldots, \bs{z}^{(n)}$ are vectors in $\R^m$ such that for $1 \le t \le n$, the two vectors $\bs{z}^{(t-1)}$ and $\bs{z}^{(t)}$ differ only in a few components, and our task is to compute all antitonic (i.e.\ monotone decreasing) approximations $A(\bs{z}^{(0)}), A(\bs{z}^{(1)}), \ldots, A(\bs{z}^{(n)})$. We show that $A(\bs{z}^{(t)})$ can be computed efficiently, provided we know already $A(\bs{z}^{(t-1)})$. Briefly speaking, this is achieved by noticing that $A(\bs{z}^{(t-1)})$ and $A(\bs{z}^{(t)})$ share some identical components, and that the remaining components of $A(\bs{z}^{(t)})$ can be determined directly from $A(\bs{z}^{(t-1)})$ and $\bs{z}^{(t)}$ with only a few operations.

The efficient computation of a sequence of antitonic approximations appears naturally in the context of isotonic distributional regression, see \cite{Henzi2021}, \cite{Moesching_Duembgen_2020} and \cite{Jordan2021}. There, one observes random pairs $(X_1,Y_1),$ $(X_2,Y_2), \ldots, (X_n,Y_n)$ in $\XX \times \R$ such that, conditional on $(X_i)_{i=1}^n$, the random variables $Y_1, Y_2, \ldots, Y_n$ are independent with distribution functions $F_{X_1}, F_{X_2}, \ldots, F_{X_n}$, where $(F_x)_{x \in \XX}$ is an unknown family of distribution functions. Then the goal is to estimate the latter family under the sole assumption that $F_{x} \ge F_{x'}$ pointwise whenever $x \preceq x'$. This notion of ordering of distributions is known as stochastic ordering, or first order stochastic dominance. This isotonic distributional regression leads to the aforementioned least squares problem, where $x_1, \ldots, x_m$ denote the different elements of $\{X_1,X_2,\ldots,X_n\}$, and $\bs{z}^{(t)}$ has components
\[
	z_j^{(t)} \
		:= \ w_j^{-1} \sum_{i \,: \, X_i = x_j} 1_{[Y_i \le Y_{(t)}]}
\]
with $w_j := \#\{i \le n : X_i = x_j\}$, $Y_{(0)}:=-\infty$ and $Y_{(t)}$ is the $t$-th order statistic of the sample $\{Y_1,Y_2,\ldots,Y_n\}$. 

Section~\ref{sec:Antitonic.LSE} provides some facts about monotone least squares which are useful for the present task. For a complete account and derivations, we refer to \cite{Barlow_1972} and \cite{Robertson_1988}. Then it is shown in Section~\ref{sec:PAVA.etc} how to turn this into an efficient computation scheme in case of a total order $\preceq$. Finally, we discuss the specific application to isotonic distributional regression, and provide numerical experiments which show that computation times of the naive approach are decreased substantially with our procedure.

\section{Some facts about antitonic least squares estimation}
\label{sec:Antitonic.LSE}

Since the sum on the right hand side of \eqref{eq:A(z)} is a strictly convex and coercive function of $\bs{f} \in \R^m$, and since $\Rmdx$ is a closed and convex set, $A(\bs{z})$ is well-defined. It possesses several well-known characterizations, two of which are particularly useful for our considerations.

The first characterization uses local weighted averages. Let us first introduce some notations. In this article, upper, lower and level sets are seen as subsets of $\{1, \ldots, m\}$ inheriting the structure of $(\XX,\preceq)$. More precisely, a set $U \subset \{1, \ldots, m\}$ is an upper set if $i \in U$ and  $x_i \preceq x_j$ imply that $j \in U$. A set $L \subset \{1, \ldots, m\}$ is a lower set if $j \in L$ and  $x_i \preceq x_j$ imply that $i \in L$. The families of all upper and all lower sets are denoted by $\UU$ and $\LL$, respectively. For a non-empty set $S\subset\{1,\ldots,m\}$, its weight and the weighted average of $\bs{z}$ over $S$ are respectively defined as
\[
	w_S^{} \
		:= \ \sum_{j\in S} w_j
	\quad\text{and}\quad
	M_S^{}(\bs{z}) \
		:= \ w_S^{-1}\sum_{j\in S} w_jz_j .
\]

\paragraph{Characterization~I.}
For any index $1 \le j \le m$,
\[
	A_j(\bs{z}) \
	= \ \min_{U \in \UU : \, j \in U} \, 
		\max_{L \in \LL : \, j \in L} \,
			M_{U \cap L}(\bs{z}) \
	= \ \max_{L \in \LL : \, j \in L} \,
		\min_{U \in \UU : \, j \in U} \,
			M_{L \cap U}(\bs{z}) .
\]

For all vectors $\bs{f}\in \mathbb{R}^m$, numbers $\xi\in\R$ and relations $\ltimes$ in $\{<,\leq,=,\geq,>\}$, let
\[
	[\bs{f} \ltimes \xi]
	\ := \
	\{j \in \{1, \ldots, m\}: f_j \ltimes \xi\}.
\]
For example, the family of sets $[\bs f = \xi]$ indexed by $\xi \in \{f_1, \ldots,f_m\}$ yields a partition of $\{1, \ldots, m\}$ such that two indices $i$ and $j$ belong to the same block if and only if $f_i = f_j$. In case of $\bs{f} \in \Rmdx$, $[\bs f < \xi]$ and $[\bs f \leq \xi]$ are upper sets, whereas $[\bs f > \xi]$ and $[\bs f \geq \xi]$ are lower sets.

\paragraph{Characterization~II.}
A vector $\bs{f} \in \Rmdx$ equals $A(\bs{z})$ if and only if for any number $\xi \in \{f_1,\ldots,f_m\}$,
\begin{align}
	M_{U \cap [\bs{f} = \xi]}(\bs{z}) \
		&\ge \ \xi
			\quad \text{for $U \in \UU$ such that} \
				U\cap [\bs{f} = \xi] \ne \emptyset , \label{Eq:CharacII_1}\\
	M_{L \cap [\bs{f} = \xi]}(\bs{z}) \
		& \le \ \xi
			\quad \text{for $L \in \LL$ such that} \
				L\cap [\bs{f} = \xi] \ne \emptyset .\label{Eq:CharacII_2}
\end{align}
In particular, $\xi = M_{[\bs{f} = \xi]}(\bs{z})$.

One possible reference for Characterizations I and II is \cite{Dominguez-Menchero_2007}. They treat the case of a quasi-order $\preceq$ and more general target functions $\sum_{j=1}^m h_j(f_j)$ to be minimized over $\bs{f} \in \Rmdx$. For the present setting with an arbitrary binary relation $\preceq$ and weighted least squares, a relatively short and self-contained derivation of these two characterizations is available from the authors upon request.

The next lemma summarizes some facts about changes in $A(\bs{z})$ if some components of $\bs{z}$ are increased.

\begin{Lemma}
\label{lem:partial.orders}
Let $\bs{z}, \bs{\tilde{z}} \in \R^m$ such that $\bs{\tilde{z}} \ge \bs{z}$ component-wise. Then the following conclusions hold true for $\bs f := A(\bs z)$, $\bs{\tilde{f}}:=A(\bs{\tilde{z}})$ and $K := \{k : \tilde{z}_k > z_k\}$:
\begin{description}
\item[(i)] $\bs{f} \le \bs{\tilde{f}}$ component-wise.
\item[(ii)] $\tilde{f}_i = f_i$ whenever $f_i < \min_{k \in K} f_k$.
\item[(iii)] $\tilde{f}_i = f_i$ whenever $\tilde{f}_i > \max_{k \in K} \tilde{f}_k$.
\item[(iv)] $\tilde{f}_i = \tilde{f}_j$ whenever $f_i = f_j$ and $x_i, x_j \preceq x_k$ for all $k \in K$.
\end{description}
\end{Lemma}

\begin{figure}
\centering
\includegraphics[width = 0.8\textwidth]{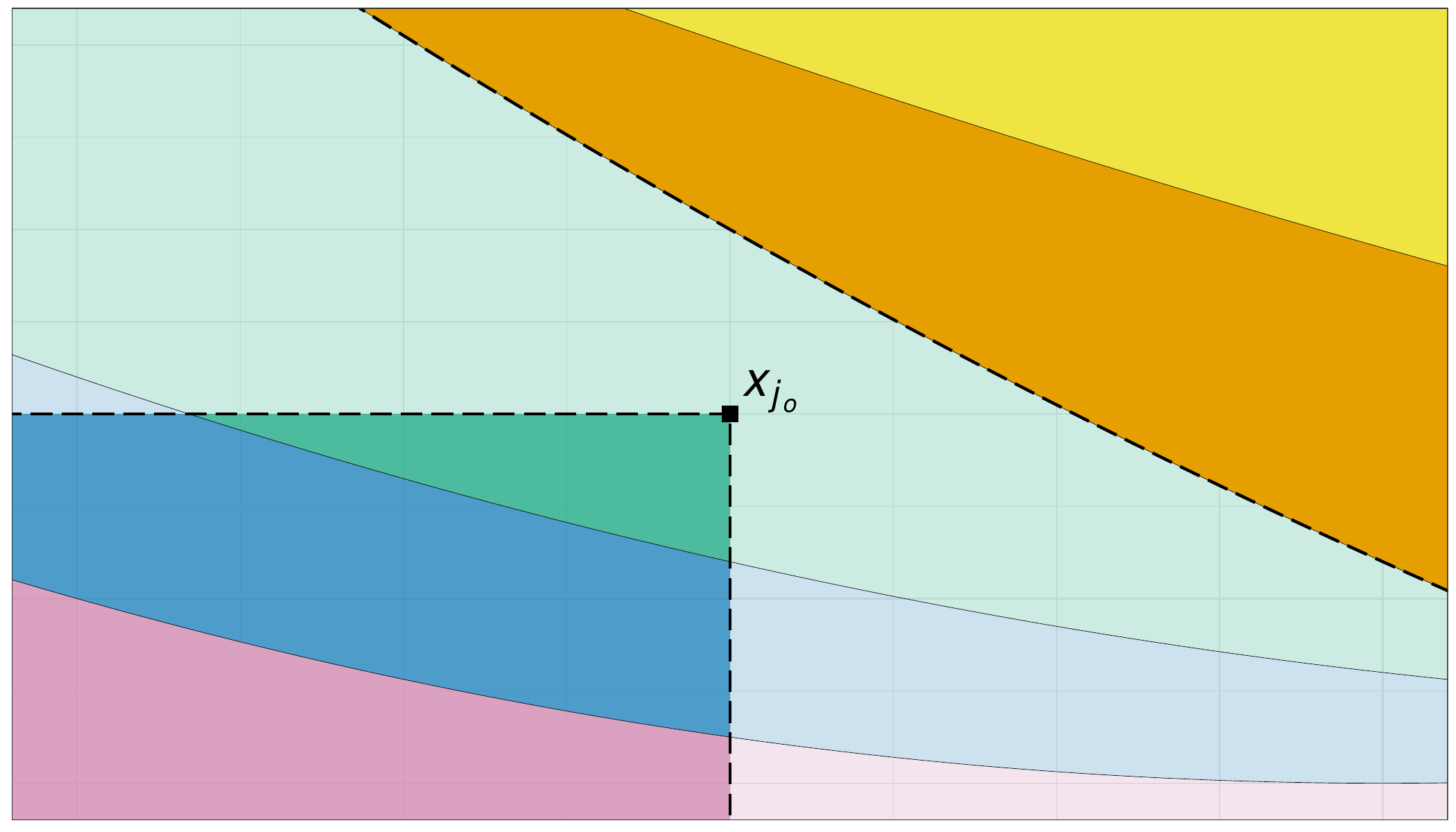} 
\caption{Illustration of the statements of Lemma~\ref{lem:partial.orders} on $\mathbb{R}^2$.}
\label{fig:illustration}
\end{figure}

Figure~\ref{fig:illustration} illustrates the statements of Lemma~\ref{lem:partial.orders} on $\mathbb{R}^2$ equipped with the componentwise order in case of $K = \{j_o\}$. The colored areas show level sets of a hypothetical antitonic regression $\bs{f}$, and $x_{j_o}$ is the point where $\tilde{z}_{j_o} > z_{j_o}$. By part~(ii) of Lemma~\ref{lem:partial.orders}, we know that $\tilde{f}_i = f_i$ if $f_i < f_{j_o}$, so the values of $\bs{f}$ and $\bs{\tilde{f}}$ are equal on the orange and yellow regions in the top right corner, which is indicated by saturated colors. Furthermore, when passing from $\bs{z}$ to $\bs{\tilde{z}}$, the slightly transparent pink, blue and green level sets on the bottom left (including the point $x_{j_o}$) can only be merged, but never be split. This follows from part~(iv) of Lemma~\ref{lem:partial.orders}. Finally, for all points in the faded pink, blue and green areas, there is no statement about the behavior of the antitonic regression when passing from $\bs{z}$ to $\bs{\tilde{z}}$.

\pagebreak

\begin{proof}[\bf{Proof of Lemma~\ref{lem:partial.orders}}]
Part~(i) is a direct consequence of Characterization~I.

As to part~(ii), if $f_i < \min_{k \in K} f_k$, then $K \subset [\bs{f} > f_i]$, whence
\begin{align*}
	\tilde{f}_i \
	&= \ \max_{L \in \LL : \, i \in L} \,
		 \min_{U \in \UU : \, i \in U} \, M_{L\cap U}^{}(\bs{\tilde{z}})
		 & (\text{Char.~I}) \\
	&\le \ \max_{L \in \LL : \, i \in L} \, M_{L \cap [\bs{f} \le f_i]}^{}(\bs{\tilde{z}})
		& (i \in [\bs{f} \le f_i] \in \UU) \\
	&= \ \max_{L \in \LL : \, i \in L} \, M_{L \cap [\bs{f} \le f_i]}^{}(\bs{z})
		& (K \cap [\bs{f} \le f_i] = \emptyset) \\
	&= \ \max_{L \in \LL : \, i \in L}
		\sum_{\xi \le f_i: \, L \cap [\bs{f} = \xi] \ne \emptyset}
			\frac{w_{L \cap [\bs{f} = \xi]}}{w_{L\cap [\bs{f} \le f_i]}} \,
				M_{L \cap [\bs{f} = \xi]}^{}(\bs{z}) \\
	&\le \ \max_{L \in \LL : \, i \in L}
		\sum_{\xi \le f_i: \, L \cap [\bs{f} = \xi] \ne \emptyset}
			\frac{w_{L \cap [\bs{f} = \xi]}}{w_{L\cap [\bs{f} \le f_i]}} \,
				\xi
		& (\text{Char.~II}) \\
	&\le \ f_i .
\end{align*}
This inequality and part~(i) show that $\tilde{f}_i = f_i$.

Part~(iii) is proved analogously. If $\tilde{f}_i > \max_{k \in K} \tilde{f}_k$, then $K \subset [\bs{\tilde{f}} < \tilde{f}_i]$, whence
\begin{align*}
	f_i \
	&= \ \min_{U \in \UU : \, i \in U} \,
		\max_{L \in \LL : \, i \in L} \, M_{U\cap L}^{}(\bs{z})
		 & (\text{Char.~I}) \\
	&\ge \ \min_{U \in \UU : \, i \in U} \,
		M_{U \cap [\bs{\tilde{f}} \ge \tilde{f}_i]}^{}(\bs{z})
		& (i \in [\bs{\tilde{f}} \ge \tilde{f}_i] \in \LL) \\
	&= \ \min_{U \in \UU : \, i \in U} \,
		M_{U \cap [\bs{\tilde{f}} \ge \tilde{f}_i]}^{}(\bs{\tilde{z}})
		& (K \cap [\bs{\tilde{f}} \ge \tilde{f}_i] = \emptyset) \\
	&= \ \min_{U \in \UU : \, i \in U}
		\sum_{\xi \ge \tilde{f}_i: \, U \cap [\bs{\tilde{f}} = \xi] \ne \emptyset}
			\frac{w_{U \cap [\bs{\tilde{f}} = \xi]}}
				{w_{U\cap [\bs{\tilde{f}} \ge \tilde{f}_i]}} \,
				M_{U \cap [\bs{\tilde{f}} = \xi]}^{}(\bs{\tilde{z}}) \\
	&\ge \ \min_{U \in \UU : \, i \in U}
		\sum_{\xi \ge \tilde{f}_i: \, U \cap [\bs{\tilde{f}} = \xi] \ne \emptyset}
			\frac{w_{U \cap [\bs{\tilde{f}} = \xi]}}
				{w_{U\cap [\bs{\tilde{f}} \le \tilde{f}_i]}} \,
				\xi
		& (\text{Char.~II}) \\
	&\ge \ \tilde{f}_i .
\end{align*}
This inequality and part~(i) show that $\tilde{f}_i = f_i$.

Part~(iv) follows directly from parts~(i) and (iii). Let $i$ and $j$ be different indices such that $f_i = f_j$ and $x_i, x_j \preceq x_k$ for all $k \in K$. It follows from $\bs{\tilde{f}} \in \Rmdx$ that $\tilde{f}_i, \tilde{f}_j \ge \max_{k \in K} \tilde{f}_k$. Consequently, if $\tilde{f}_j > \tilde{f}_i$, then $\tilde{f}_j > \max_{k \in K} \tilde{f}_k$, so parts~(i) and (iii) would imply that
\[
	\tilde{f}_i \ \ge \ f_i = f_j \ = \ \tilde{f}_j ,
\]
contradicting $\tilde{f}_j > \tilde{f}_i$.
\end{proof}

\paragraph{The special case of a total order.}
If one replaces the binary relation $\preceq$ by a total order $\leq$ on $\XX$, as for example in the case of the usual total order on a subset of $\R$, the conclusions of Lemma~\ref{lem:partial.orders} take a simpler form. In case of a total order, we assume that the covariates are ordered as follows
\[
	x_1 \ \le \ x_2 \ \le \ \cdots \ \le \ x_m ,
\]
so that $i\leq j$ implies that $x_i\leq x_j$, while $x_i < x_j$ implies that $i < j$. 

\begin{Corollary}
\label{cor:total.order}
Let $\bs{z}, \bs{\tilde{z}} \in \R^m$ such that $\bs{z} \le \bs{\tilde{z}}$ component-wise. Then the following conclusions hold true for $\bs f := A(\bs z)$ and $\bs{\tilde{f}}:=A(\bs{\tilde{z}})$:
\begin{description}
\item[(i)] $\bs{f} \le \bs{\tilde{f}}$ component-wise.

\item[(ii)] Let $k \in \{1,\ldots,m-1\}$ such that $f_k > f_{k+1}$ and $(\tilde{z}_j)_{j > k} = (z_j)_{j > k}$. Then
\[
	(\tilde{f}_j)_{j > k} \ = \ (f_j)_{j > k} .
\]

\item[(iii)] Let $k \in \{2,\ldots,m\}$ such that $\tilde{f}_{k-1} > \tilde{f}_k$ and $(\tilde{z}_j)_{j < k} = (z_j)_{j < k}$. Then
\[
	(\tilde{f}_j)_{j < k} \ = \ (f_j)_{j < k} .
\]

\item[(iv)] Let $k \in \{2,\ldots,m\}$ such that $(\tilde{z}_j)_{j < k} = (z_j)_{j < k}$.
Then
\[
	\{j < k : \tilde{f}_j > \tilde{f}_{j+1}\}
	\ \subset \
	\{j < k : f_j > f_{j+1}\} .
\]
\end{description}
\end{Corollary}

\section{A sequential algorithm for total orders}
\label{sec:PAVA.etc}

Lemma~\ref{lem:partial.orders} is potentially useful to accelerate algorithms for isotonic distributional regression with arbitrary partial orders, possibly in conjunction with the recursive partitioning algorithm by \citet{Luss2014}, but this will require additional research. Now we focus on improvements of the well-known pool-adjacent-violators algorithm (PAVA) for a total order.

\subsection{General considerations}

In what follows, we assume that $x_1 < \cdots < x_m$, so $\Rmdx$ coincides with $\Rmda = \bigl\{ \bs{f} \in \R^m : f_1 \ge \cdots \ge f_m\}$. To understand the different variants of the PAVA, let us recall two basic facts about $A(\bs{z})$. Let $\PP = (P_1,\ldots,P_d)$ be a partition of $\{1,\ldots,m\}$ into blocks $P_s = \{b_{s-1}+1,\ldots,b_s\}$, where $0 = b_0 < b_1 < \cdots < b_d = m$, and let $\R^m_{\PP}$ be the set of vectors $\bs{f} \in \R^m$ such that $f_i = f_j$ whenever $i,j$ belong to the same block of $\PP$.

\noindent
\textbf{Fact~1.} Let $r_1 > \cdots > r_d$ be the sorted elements of $\{A_i(\bs{z}) : 1 \le i \le m\}$, and let $\PP$ consist of the blocks $P_s = \{i : A_i(\bs{z}) = r_s\}$. Then $r_s = M_{P_s}^{}(\bs{z})$ for $1 \le s \le d$.

\noindent
\textbf{Fact~2.} Suppose that $A(\bs{z}) \in \R^m_{\PP}$ for a given partition $\PP$ with $d \ge 2$ blocks. If $s \in \{1,\ldots,d-1\}$ such that $M_{P_s}^{}(\bs{z}) \le M_{P_{s+1}}^{}(\bs{z})$, then $A_i(\bs{z})$ is constant in $i \in P_s \cup P_{s+1}$. That means, one may replace $\PP$ with a coarser partition by pooling $P_s$ and $P_{s+1}$ and still, $A(\bs{z}) \in \R^m_{\PP}$.

Fact~1 is a direct consequence of Characterization~II. To verify Fact~2, suppose that $\bs{f} \in \Rmda \cap \R^m_{\PP}$ such that $f_i = r_s$ for $i \in P_s$, $f_i = r_{s+1}$ for $i \in P_{s+1}$, and $r_s > r_{s+1}$. Now we show that $\bs f$ cannot be equal to $A(\bs z)$. For $t \ge 0$ let $\bs{f}(t) \in \R^m_{\PP}$ be given by
\[
	f_i(t) \ = \ f_i
		- 1_{[i \in P_s]} t w_{P_s}^{-1}
		+ 1_{[i \in P_{s+1}]} t w_{P_{s+1}}^{-1} .
\]
Then $\bs{f}(0) = \bs{f}$, and $\bs{f}(t) \in \Rmda$ if $t \le (r_s - r_{s+1})/(w_{P_{s+1}}^{-1} + w_{P_s}^{-1})$. But
\[
	\frac{d}{dt} \Big|_{t = 0} \sum_{i=1}^m w_i(f_i(t) - z_i)^2
	\ = \ 2 (r_{s+1} - r_s) - 2 \bigl( M_{P_{s+1}}^{}(\bs{z}) - M_{P_s}^{}(\bs{z}) \bigr)
	\ < \ 0 ,
\]
so for sufficiently small $t>0$, $\bs{f}(t)\in \Rmda$ and is superior to $\bs{f}(0)$. Hence $\bs{f} \ne A(\bs{z})$.

Facts~1 and 2 indicate already a general PAV strategy to compute $A(\bs{z})$. One starts with the finest partition $\PP = (\{1\}, \ldots, \{m\})$. As long as $\PP$ contains two neighboring blocks $P_s$ and $P_{s+1}$ such that $M_{P_s}^{}(\bs{z}) \ge M_{P_{s+1}}^{}(\bs{z})$, the partition $\PP$ is coarsened by replacing $P_s$ and $P_{s+1}$ with the block $P_s \cup P_{s+1}$.

\paragraph{Standard PAVA.}
Specifically, one works with three tuples: $\PP = (P_1,\ldots,P_d)$ is a partition of $\{1,\ldots,b_d\}$ into blocks $P_s = \{b_{s-1}+1,\ldots,b_s\}$, where $0 = b_0 < b_1 < \cdots < b_d$. The number $b_d$ is running from $1$ to $m$, and the number $d \ge 1$ changes during the algorithm, too. The tuples $\WW = (W_1,\ldots,W_d)$ and $\MM = (M_1,\ldots,M_d)$ contain the corresponding weights $W_s = w_{P_s}^{}$ and weighted means $M_s = M_{P_s}^{}(\bs{z})$. Before increasing $b_d$, the tuples $\PP$, $\WW$ and $\MM$ describe the minimizer of $\sum_{i=1}^{b_d} w_i (f_i - z_i)^2$ over all $\bs{f} \in \R^{b_d}_{\downarrow}$. Here is the complete algorithm:

\noindent
Initialization: We set $\PP \leftarrow (\{1\})$, $\WW \leftarrow (w_1)$, $\MM \leftarrow (z_1)$, and $d \leftarrow 1$.

\noindent
Induction step: If $b_d < m$, we add a new block by setting
\[
	\PP \ \leftarrow \ (\PP,\{b_d+1\}) ,
	\quad
	\WW \ \leftarrow \ (\WW,w_{b_d+1}^{}) ,
	\quad
	\MM \ \leftarrow \ (\MM, z_{b_d+1}^{}) ,
\]
and $d \leftarrow d+1$. Then, while $d > 1$ and $M_{d-1} \le M_d$, we pool the ``violators'' $P_{d-1}$ and $P_d$ by setting
\begin{align*}
	\PP \ &\leftarrow \ \bigl( (P_j)_{j < d-1},
		P_{d-1} \cup P_d \bigr) , \\
	\MM \ &\leftarrow \ \Bigl( (W_j)_{j < d-1},
		\frac{W_{d-1}M_{d-1} + W_dM_d}{W_{d-1} + W_d} \Bigr) , \\
	\WW \ &\leftarrow \ \bigl( (W_j)_{j < d-1},
		W_{d-1} + W_d \bigr) ,
\end{align*}
and $d \leftarrow d-1$.

\noindent Finalization: Eventually, $\PP$ is a partition of $\{1,\ldots,m\}$ into blocks such that $M_1 > \cdots > M_d$ and
\[
	A_j(\bs{z}) \ = \ M_s \quad\text{for} \ j \in P_s \ \text{and}\ 1 \le s \le d .
\]

\paragraph{Modified PAVA.}
In our specific applications of the PAVA, we are dealing with vectors $\bs{z}$ containing larger blocks $\{a,\ldots,b\}$ on which $i \mapsto z_i$ is constant. Indeed, in regression settings with continuously distributed covariates and responses, $\bs{z}$ will always be a $\{0,1\}$-valued vector. Then it is worthwhile to utilize fact~2 and modify the initialization as well as the very beginning of the induction step as follows:

For the initialization, we determine the largest index $b_1$ such that $z_1 = \cdots = z_{b_1}$ and the corresponding weight $W_{P_1}$ with $P_1 = \{1,\ldots,b_1\}$. Then we set $\PP \leftarrow (P_1)$, $\WW \leftarrow (w_{P_1}^{})$ and $\MM \leftarrow (z_{b_1}^{})$, where $P_1 = \{1,\ldots,b_1\}$.

At the beginning of the induction step, we determine the largest index $b_{d+1} > b_d$ such that $z_{b_d+1} = \cdots = z_{b_{d+1}}$ and the corresponding weight $W_{P_{d+1}}^{}$ with $P_{d+1} = \{b_d+1,\ldots,b_{d+1}\}$. Then we set $\PP \leftarrow (\PP,P_{d+1})$, $\WW \leftarrow (\WW,W_{P_{d+1}})$, $\MM \leftarrow (\MM, z_{b_{d+1}}^{})$, and $d \leftarrow d+1$.

\paragraph{Abridged PAVA.}
Suppose that we have computed $A(\bs{z})$ with corresponding tuples $\PP = (P_1,\ldots,P_d)$, $\WW = (W_1,\ldots,W_d)$ and $\MM = (M_1,\ldots,M_d)$ via the PAVA. Now let $\bs{\tilde{z}} \in \R^m$ such that $\tilde{z}_{j_o} > z_{j_o}$ for one index $j_o \in \{1,\ldots,m\}$, while $(\tilde{z}_j)_{j \ne j_o} = (z_j)_{j \ne j_o}$. Let $j_o \in P_{s_o}$ with $s_o \in \{1,\ldots,d\}$. By parts~(ii) and (iv) of Corollary~\ref{cor:total.order}, the partition corresponding to $A(\bs{\tilde{z}})$ will be a coarsening of the partition with the following blocks:
\[
	P_s \ \ \text{for} \ 1 \le s < s_o ,
	\quad
	\{b_{s_o-1} + 1,\ldots,j_o\} ,
	\quad
	\{j\} \ \ \text{for} \ j_o < j \le b_{s_o} ,
	\quad
	P_s \ \ \text{for} \ s_o < s \le d .
\]
Moreover, $A_i(\bs{\tilde{z}}) = A_i(\bs{z})$ for $i > b_{s_o}$. This allows us to compute $A(\bs{\tilde{z}})$ as follows, keeping copies of the auxiliary objects for $A(\bs{z})$ and indicating this with a superscript $\bs{z}$:

\noindent
Initialization: We determine $s_o \in \{1,\ldots,d^{\bs{z}}\}$ such that $j_o \in P_{s_o}^{\bs{z}}$. Then we set
\begin{align*}
	\PP \ &\leftarrow \ \bigl( (P_s^{\bs{z}})_{s < s_o}^{},
		\{b_{s_o-1}^{\bs{z}}+1,\ldots,j_o\} \bigr) , \\
	\MM \ &\leftarrow \ \bigl( (M_s^{\bs{z}})_{s < s_o}^{},
		M_{P_{s_o}}^{}(\bs{\tilde{z}}) \bigr) , \\
	\WW \ &\leftarrow \ \bigl( (W_s^{\bs{z}})_{s < s_o}^{},
		w_{P_{s_o}}^{} \bigr)
\end{align*}
and $d \leftarrow s_o$. While $d > 1$ and $M_{d-1} \le M_d$, we pool the violators $P_{d-1}$ and $P_d$ as in the induction step of PAVA. (This initialization is justified by part~(iv) of Corollary~\ref{cor:total.order}.)

\noindent
Induction step: If $j_o < b_{s_o}^{\bs{z}}$, we run the induction step of PAVA for $b_d$ running from $j_o+1$ to $b_{s_o}^{\bs{z}}$ with $\bs{\tilde{z}}$ in place of $\bs{z}$.

\noindent
Finalization: If $b_{s_o}^{\bs{z}} < m$, we set
\begin{align*}
	\PP \ &\leftarrow \
		\bigl( \PP, (P_s^{\bs{z}})_{s_o < s \le d^{\bs{z}}} \bigr) , \\
	\MM \ &\leftarrow \
		\bigl( \MM, (M_s^{\bs{z}})_{s_o < s \le d^{\bs{z}}} \bigr) , \\
	\WW \ &\leftarrow \
		\bigl( \WW, (W_s^{\bs{z}})_{s_o < s \le d^{\bs{z}}} \bigr)
\end{align*}
and $d \leftarrow d + d^{\bs{z}} - s_o$. The new pair $(\PP,\MM)$ yields the vector $A(\bs{\tilde{z}})$. This finalization is justified by part~(ii) of Corollary~\ref{cor:total.order}.

\paragraph{Computational complexity.}
It directly follows from the algorithmic description that when $A(\bs{z})$ is available, the abridged PAVA for computing $A(\bs{\tilde{z}})$ requires not more operations than the standard PAVA. Its computational complexity is therefore at most of order $O(m)$ if $x_1, \dots, x_m$ are already sorted. More precisely, the number of averaging operations in the abridged PAVA is bounded from above by $d^{\bs{z}} + (b_{s_o}^{\bs{z}} - b_{s_o - 1}^{\bs{z}})$, where $d^{\bs{z}}$ is the partition size of the antitonic regression $A(\bs{z})$ and $b_{s_o}^{\bs{z}} - b_{s_o - 1}^{\bs{z}}$ is the number of elements in the set $P_{s_o}^{\bs{z}}$ containing the index $j_o$ where the value of $\bs{z}$ changes. In many practical applications this number is much smaller than $m$, but in the worst case it may equal exactly $m$; for example, let $w_i = 1$ and $z_i = m - i$ for $i = 1, \dots, m$, $j_o = m$, and $\tilde{z}_m = m^2$.

\paragraph{Numerical example.}
We illustrate the previous procedures with two vectors $\bs{z}, \bs{\tilde{z}} \in \R^9$ and $\bs{w} = (1)_{j=1}^9$. Table~\ref{tab:PAVA1} shows the main steps of the PAVA for $\bs{z}$. The first line shows the components of $\bs{z}$, the other lines contain the current candidate for $(f_j)_{j=1}^{b_d}$, where $\bs{f} = A(\bs{z})$ eventually, and the current partition $\PP$ is indicated by extra vertical bars. Table~\ref{tab:PAVA2} shows the abridged PAVA for two different vectors $\bs{\tilde{z}}$.

\begin{table}
\[
	\begin{array}{|c||c|c|c|c|c|c|c|c|c||l|}
	\cline{1-10}
	\bs{z} &\ 1 &\ 3 &\ 2 &  0 & -1 &  1 & 1/2& -1 &\ 1 &
		\multicolumn{1}{c}{} \\
	\hline\hline\hline
	b_d=1   &  1 & \multicolumn{8}{c||}{} &
	d = 1 \\
	\hline\hline
	b_d=2   &  1 &  3 & \multicolumn{7}{c||}{} &
	d = 2 \\
	\cline{2-11}
	       & \multicolumn{1}{|c}{2}
	            & \multicolumn{1}{c|}{2}
	                 & \multicolumn{7}{c||}{} &
	d = 1 \\
	\hline\hline
	b_d=3   & \multicolumn{1}{|c}{2}
	            & \multicolumn{1}{c|}{2}
	                 &  2
	                      & \multicolumn{6}{c||}{} &
	d = 2 \\
	\cline{2-11}
	       & \multicolumn{1}{|c}{2}
	            & \multicolumn{1}{c}{2}
	                 & \multicolumn{1}{c|}{2}
	                      & \multicolumn{6}{c||}{} &
	d = 1 \\
	\hline\hline
	b_d = 4 & \multicolumn{1}{|c}{2}
	            & \multicolumn{1}{c}{2}
	                 & \multicolumn{1}{c|}{2}
	                      &  0
	                           & \multicolumn{5}{c||}{} &
	d = 2 \\
	\hline\hline
	b_d = 5 & \multicolumn{1}{|c}{2}
	            & \multicolumn{1}{c}{2}
	                 & \multicolumn{1}{c|}{2}
	                      &  0 & -1
	                                & \multicolumn{4}{c||}{} &
	d = 3 \\
	\hline\hline
	b_d = 6 & \multicolumn{1}{|c}{2}
	            & \multicolumn{1}{c}{2}
	                 & \multicolumn{1}{c|}{2}
	                      &  0 & -1 &  1
	                                     & \multicolumn{3}{c||}{} &
	d = 4 \\
	\cline{2-11}
	       & \multicolumn{1}{|c}{2}
	            & \multicolumn{1}{c}{2}
	                 & \multicolumn{1}{c|}{2}
	                      &  0 & \multicolumn{1}{|c}{0}
	                                & \multicolumn{1}{c|}{0}
	                                     & \multicolumn{3}{c||}{} &
	d = 3 \\
	\cline{2-11}
	       & \multicolumn{1}{|c}{2}
	            & \multicolumn{1}{c}{2}
	                 & \multicolumn{1}{c|}{2}
	                      & \multicolumn{1}{|c}{0}
	                           & \multicolumn{1}{c}{0}
	                                & \multicolumn{1}{c|}{0}
	                                     & \multicolumn{3}{c||}{} &
	d = 2 \\
	\hline\hline
	b_d = 7 & \multicolumn{1}{|c}{2}
	             & \multicolumn{1}{c}{2}
	                   & \multicolumn{1}{c|}{2}
	                         & \multicolumn{1}{|c}{0}
	                               & \multicolumn{1}{c}{0}
	                                     & \multicolumn{1}{c|}{0}
	                                           & 1/2
	                                                 & \multicolumn{2}{c||}{} &
	d = 3 \\
	\cline{2-11}
	       & \multicolumn{1}{|c}{2}
	            & \multicolumn{1}{c}{2}
	                 & \multicolumn{1}{c|}{2}
	                      & \multicolumn{1}{|c}{1/8}
	                           & \multicolumn{1}{c}{1/8}
	                                & \multicolumn{1}{c}{1/8}
	                                     & \multicolumn{1}{c|}{1/8}
	                                          & \multicolumn{2}{c||}{} &
	d = 2 \\
	\hline\hline
	b_d = 8 & \multicolumn{1}{|c}{2}
	            & \multicolumn{1}{c}{2}
	                 & \multicolumn{1}{c|}{2}
	                      & \multicolumn{1}{|c}{1/8}
	                           & \multicolumn{1}{c}{1/8}
	                                & \multicolumn{1}{c}{1/8}
	                                     & \multicolumn{1}{c|}{1/8}
	                                          & -1 &    &
	d = 3 \\
	\hline\hline
	b_d = 9 & \multicolumn{1}{|c}{2}
	            & \multicolumn{1}{c}{2}
	                 & \multicolumn{1}{c|}{2}
	                      & \multicolumn{1}{|c}{1/8}
	                           & \multicolumn{1}{c}{1/8}
	                                & \multicolumn{1}{c}{1/8}
	                                     & \multicolumn{1}{c|}{1/8}
	                                          & -1 &  1 &
	d = 4 \\
	\cline{2-11}
	       & \multicolumn{1}{|c}{2}
	            & \multicolumn{1}{c}{2}
	                 & \multicolumn{1}{c|}{2}
	                      & \multicolumn{1}{|c}{1/8}
	                           & \multicolumn{1}{c}{1/8}
	                                & \multicolumn{1}{c}{1/8}
	                                     & \multicolumn{1}{c|}{1/8}
	                                          & \multicolumn{1}{|c}{0}
	                                               & \multicolumn{1}{c||}{0} &
	d = 3 \\
	\hline
	\end{array}
\]
\caption{Running the PAVA for a vector $\bs{z}$.}
\label{tab:PAVA1}
\end{table}

\begin{table}
\[
	\begin{array}{|c||c|c|c|c|c|c|c|c|c||l|}
	\cline{1-10}
	\bs{z} &\ 1 &\ 3 &\ 2 &  0 & -1 &  1 & 1/2& -1 &\ 1 &
		\multicolumn{1}{c}{} \\
	\cline{1-10}
	A(\bs{z})
	       & \multicolumn{1}{|c}{\ 2}
	            & \multicolumn{1}{c}{\ 2}
	                 & \multicolumn{1}{c|}{\ 2}
	                      & \multicolumn{1}{|c}{1/8}
	                           & \multicolumn{1}{c}{1/8}
	                                & \multicolumn{1}{c}{1/8}
	                                     & \multicolumn{1}{c|}{1/8}
	                                          & \multicolumn{1}{|c}{\ 0}
	                                               & \multicolumn{1}{c||}{\ 0} &
		\multicolumn{1}{c}{} \\
	\cline{1-10}
	\bs{\tilde{z}}
	       &\ 1 &\ 3 &\ 2 &  0 & \bs{1}
	                                &  1 & 1/2& -1 &\ 1 &
		\multicolumn{1}{c}{} \\
	\hline\hline\hline
	b_d = 5 & \multicolumn{1}{|c}{2}
	            & \multicolumn{1}{c}{2}
	                 & \multicolumn{1}{c|}{2}
	                      & \multicolumn{1}{|c}{1/2}
	                           & \multicolumn{1}{c|}{1/2}
	                                & \multicolumn{4}{c||}{} &
	d = 2 \\
	\hline\hline
	b_d = 6 & \multicolumn{1}{|c}{2}
	            & \multicolumn{1}{c}{2}
	                 & \multicolumn{1}{c|}{2}
	                      & \multicolumn{1}{|c}{1/2}
	                           & \multicolumn{1}{c|}{1/2}
	                                &  1
	                                     & \multicolumn{3}{c||}{} &
	d = 3 \\
	\cline{2-11}
	       & \multicolumn{1}{|c}{2}
	            & \multicolumn{1}{c}{2}
	                 & \multicolumn{1}{c|}{2}
	                      & \multicolumn{1}{|c}{2/3}
	                           & \multicolumn{1}{c}{2/3}
	                                & \multicolumn{1}{c|}{2/3}
	                                     & \multicolumn{3}{c||}{} &
	d = 2 \\
	\hline\hline
	b_d = 7 & \multicolumn{1}{|c}{2}
	             & \multicolumn{1}{c}{2}
	                   & \multicolumn{1}{c|}{2}
	                         & \multicolumn{1}{|c}{2/3}
	                               & \multicolumn{1}{c}{2/3}
	                                     & \multicolumn{1}{c|}{2/3}
	                                           & 1/2
	                                                 & \multicolumn{2}{c||}{} &
	d = 3 \\
	\hline\hline
	b_d = 9 & \multicolumn{1}{|c}{2}
	             & \multicolumn{1}{c}{2}
	                   & \multicolumn{1}{c|}{2}
	                         & \multicolumn{1}{|c}{2/3}
	                               & \multicolumn{1}{c}{2/3}
	                                     & \multicolumn{1}{c|}{2/3}
	                                          & 1/2
	                                               & \multicolumn{1}{|c}{0}
	                                               & \multicolumn{1}{c||}{0} &
	d = 4 \\
	\hline
	\multicolumn{11}{c}{} \\
	\multicolumn{11}{c}{} \\
	\cline{1-10}
	\bs{z} &\ 1 &\ 3 &\ 2 &  0 & -1 &  1 & 1/2& -1 &\ 1 &
		\multicolumn{1}{c}{} \\
	\cline{1-10}
	A(\bs{z})
	       & \multicolumn{1}{|c}{\ 2}
	            & \multicolumn{1}{c}{\ 2}
	                 & \multicolumn{1}{c|}{\ 2}
	                      & \multicolumn{1}{|c}{1/8}
	                           & \multicolumn{1}{c}{1/8}
	                                & \multicolumn{1}{c}{1/8}
	                                     & \multicolumn{1}{c|}{1/8}
	                                          & \multicolumn{1}{|c}{\ 0}
	                                               & \multicolumn{1}{c||}{\ 0} &
		\multicolumn{1}{c}{} \\
	\cline{1-10}
	\bs{\tilde{z}}
	       &\ 1 &\ 3 &\ 2 & \bs{2}
	                           & -1
	                                &  1 & 1/2& -1 &\ 1 &
		\multicolumn{1}{c}{} \\
	\hline\hline\hline
	b_d = 4 & \multicolumn{1}{|c}{2}
	            & \multicolumn{1}{c}{2}
	                 & \multicolumn{1}{c|}{2}
	                      &  2
	                           & \multicolumn{5}{c||}{} &
	d = 2 \\
	\cline{2-11}
	       & \multicolumn{1}{|c}{2}
	            & \multicolumn{1}{c}{2}
	                 & \multicolumn{1}{c}{2}
	                      & \multicolumn{1}{c|}{2}
	                           & \multicolumn{5}{c||}{} &
	d = 1 \\
	\hline\hline
	b_d = 5 & \multicolumn{1}{|c}{2}
	            & \multicolumn{1}{c}{2}
	                 & \multicolumn{1}{c}{2}
	                      & \multicolumn{1}{c|}{2}
	                           & -1
	                                & \multicolumn{4}{c||}{} &
	d = 2 \\
	\hline\hline
	b_d = 6 & \multicolumn{1}{|c}{2}
	            & \multicolumn{1}{c}{2}
	                 & \multicolumn{1}{c}{2}
	                      & \multicolumn{1}{c|}{2}
	                           & -1
	                                &  1
	                                     & \multicolumn{3}{c||}{} &
	d = 3 \\
	\cline{2-11}
	       & \multicolumn{1}{|c}{2}
	            & \multicolumn{1}{c}{2}
	                 & \multicolumn{1}{c}{2}
	                      & \multicolumn{1}{c|}{2}
	                           & \multicolumn{1}{|c}{0}
	                                & \multicolumn{1}{c|}{0}
	                                     & \multicolumn{3}{c||}{} &
	d = 2 \\
	\hline\hline
	b_d = 7 & \multicolumn{1}{|c}{2}
	            & \multicolumn{1}{c}{2}
	                 & \multicolumn{1}{c}{2}
	                      & \multicolumn{1}{c|}{2}
	                           & \multicolumn{1}{|c}{0}
	                                & \multicolumn{1}{c|}{0}
	                                     & 1/2
	                                          & \multicolumn{2}{c||}{} &
	d = 3 \\
	\cline{2-11}
	       & \multicolumn{1}{|c}{2}
	            & \multicolumn{1}{c}{2}
	                 & \multicolumn{1}{c}{2}
	                      & \multicolumn{1}{c|}{2}
	                           & \multicolumn{1}{|c}{1/6}
	                                & \multicolumn{1}{c}{1/6}
	                                     & \multicolumn{1}{c|}{1/6}
	                                          & \multicolumn{2}{c||}{} &
	d = 2 \\
	\hline\hline
	b_d = 9 & \multicolumn{1}{|c}{2}
	            & \multicolumn{1}{c}{2}
	                 & \multicolumn{1}{c}{2}
	                      & \multicolumn{1}{c|}{2}
	                           & \multicolumn{1}{|c}{1/6}
	                                & \multicolumn{1}{c}{1/6}
	                                     & \multicolumn{1}{c|}{1/6}
	                                               & \multicolumn{1}{|c}{0}
	                                               & \multicolumn{1}{c||}{0} &
	d = 3 \\
	\hline
	\end{array}
\]
\caption{Running the abridged PAVA for two vectors $\bs{\tilde{z}} \approx \bs{z}$.}
\label{tab:PAVA2}
\end{table}

\subsection{Application to isotonic distributional regression}

Now we consider a regression framework similar to the one discussed in \cite{Moesching_Duembgen_2020}, \cite{Henzi2021} and  \cite{Jordan2021}. We observe pairs $(X_1,Y_1),$ $(X_2,Y_2),\ldots,(X_n,Y_n)$ consisting of numbers $X_i \in \XX$ (covariate) and $Y_i \in \R$ (response), where $\XX$ is a given real interval. Conditional on $(X_i)_{i=1}^n$, the observations $Y_1, Y_2, \ldots, Y_n$ are viewed as independent random variables such that for $x \in \XX$ and $y \in \R$,
\[
	\Pr(Y_i \le y) \ = \ F_x(y)
	\quad\text{if} \ X_i = x .
\]
Here $(F_x)_{x \in \XX}$ is an unknown family of distribution functions. We only assume that $F_x(y)$ is non-increasing in $x \in \XX$ for any fixed $y \in \R$. That means, the family $(F_x)_{x \in \XX}$ is increasing with respect to stochastic order.

Let $x_1 < x_2 < \cdots < x_m$ be the elements of $\{X_1,X_2,\ldots,X_n\}$, and let
\[
	w_j \ := \ \#\{i: \, X_i = x_j\} ,
	\quad 1 \le j \le m .
\]
Then one can estimate $\bs{F}(y) := (F_{x_j}(y))_{j=1}^m$ by
\[
	\hat{\bs{F}}(y) \ := \ A(\bs{z}(y)) ,
\]
where $\bs{z}(y)$ has components
\[
	z_j(y) \ := \ w_j^{-1} \sum_{i: \, X_i = x_j} 1_{[Y_i \le y]} ,
	\quad 1 \le j \le m .
\]

Suppose we have rearranged the observations such that $Y_1 \le Y_2 \le \cdots \le Y_n$. Let $\bs{z}^{(0)} := \bs{0}$ and
\[
	\bs{z}^{(t)}
	\ := \ \Bigl( w_j^{-1} \sum_{i \le t: \, X_i = x_j} 1_{[Y_i \le Y_t]} \Bigr)_{j=1}^m
\]
for $1 \le t \le n$. Note that $\bs{z}^{(t-1)}$ and $\bs{z}^{(t)}$ differ in precisely one component, and that
\[
	\bs{z}(y) \ = \ \begin{cases}
		\bs{z}^{(0)} & \text{if} \ y < Y_1 , \\
		\bs{z}^{(t)} & \text{if} \ Y_t \le y < Y_{t+1}, \ 1 \le t < n , \\
		\bs{z}^{(n)} & \text{if} \ y \ge Y_n .
	\end{cases}
\]
Thus it suffices to compute $A(\bs{z}^{(t)})$ for $t = 0,1,\ldots,n$. But $A(\bs{z}^{(0)}) = \bs{0}$, $A(\bs{z}^{(n)}) = \bs{1}$, and for $1 \le t < n$, one may apply the abridged PAVA to the vectors $\bs{z} := \bs{z}^{(t-1)}$ and $\bs{\tilde{z}} := \bs{z}^{(t)}$. This leads to an efficient algorithm to compute all vectors $A(\bs{z}^{(t)})$, $0 \le t \le n$, if implemented properly.

\begin{figure}
\centering
\includegraphics[width=0.7\textwidth]{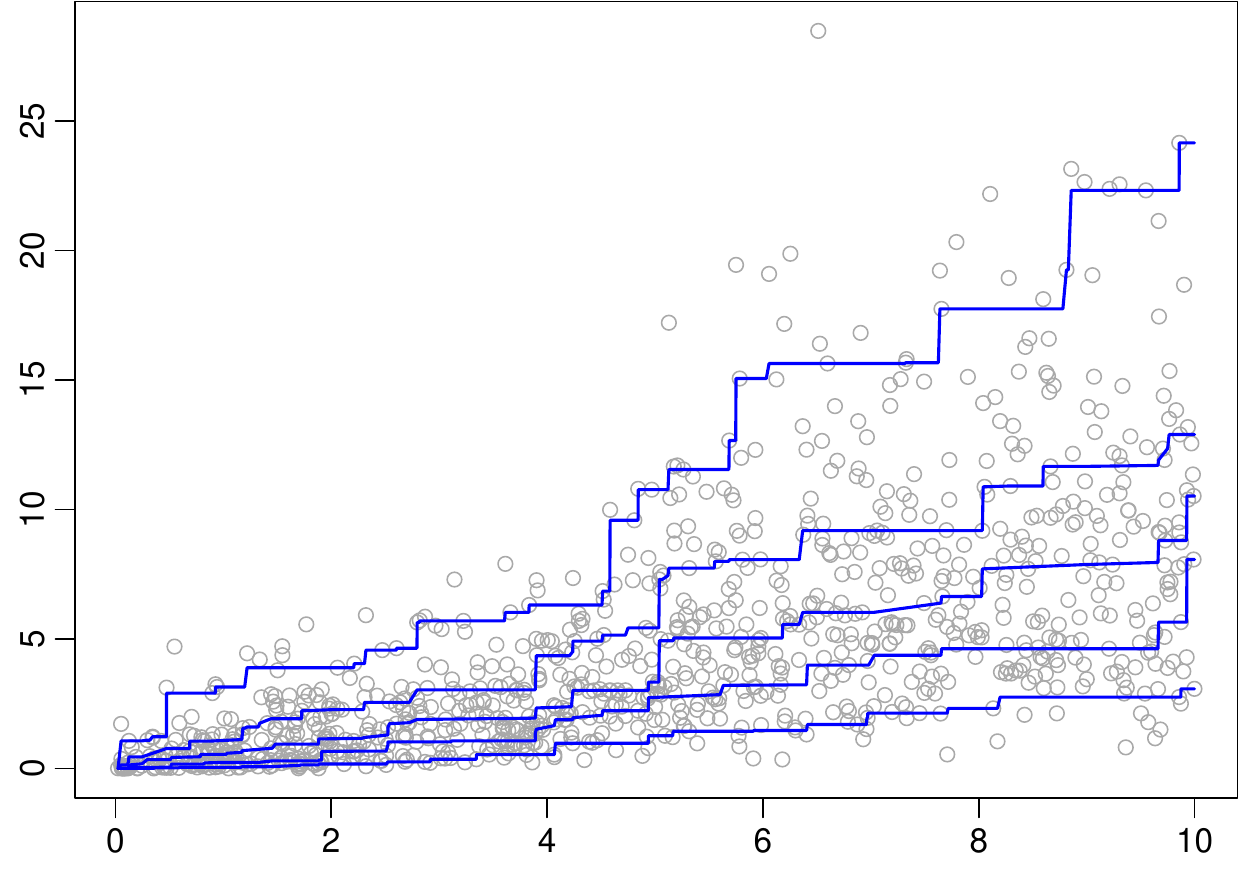}
\caption{A data set with estimated quantile curves.}
\label{fig:QuantilesYX}
\end{figure}

\paragraph{Numerical experiment 1.}
We generated data sets with $n = 1000$ independent observation pairs $(X_i,Y_i)$, $1 \le i \le n$, where $X_i$ is uniformly distributed on $[0,10]$ while $\LL(Y_i \,|\, X_i = x)$ is the gamma distribution with shape parameter $\sqrt{x}$ and scale parameter $2 + (x - 5)/\sqrt{2 + (x-5)^2}$. Figure~\ref{fig:QuantilesYX} shows one such data set. In addition, one sees estimated $\beta$-quantile curves for levels $\beta \in \{0.1, 0.25, 0.5, 0.75, 0.9\}$, resulting from the estimator $\hat{\bs{F}}$.

\begin{table}
\[
	\begin{array}{|l||rc||rc||rc|}
	\cline{1-3}
	\text{Variant of PAVA} &
	\multicolumn{2}{|l||}{\text{mean (sd) of} \ T_j} \\
	\cline{1-5}
	\text{Standard} \ \ T_1 & 6.0394 & (1.5257) &
		\multicolumn{2}{|l||}{\text{mean (sd) of} \ T_1/T_j}\\
	\hline
	\text{Modified} \ \ T_2 & 1.7482 & (0.4224) &  3.4618 & (0.3816) &
		\multicolumn{2}{|l|}{\text{mean (sd) of} \ T_2/T_3} \\
	\hline
	\text{Abridged} \ \ T_3 & 0.2080 & (0.1052) & 30.8308 & (6.1209) & 8.9012 & (1.4469) \\
	\hline
	\end{array}
\]
\caption{Computation times in seconds and ratios of running times.}
\label{tab:QuantilesYX}
\end{table}

Now we simulated $1000$ such data sets and measured the times $T_1, T_2, T_3$ for computing the estimator $\hat{\bs{F}}$ via the standard, the modified and the abridged PAVA, respectively. Table~\ref{tab:QuantilesYX} reports the sample means and standard deviations of these computation times in the $1000$ simulations. In addition, one sees the averages and standard deviations of the ratios $T_i/T_j$, for $1 \le i < j \le 3$. It turned out that using the modified instead of the standard PAVA reduced the computation time by a factor of $3.46$ already. Using the abridged PAVA yielded a further improvement by a factor of $8.90$. 

\begin{figure}
\centering
\includegraphics[width=1\textwidth]{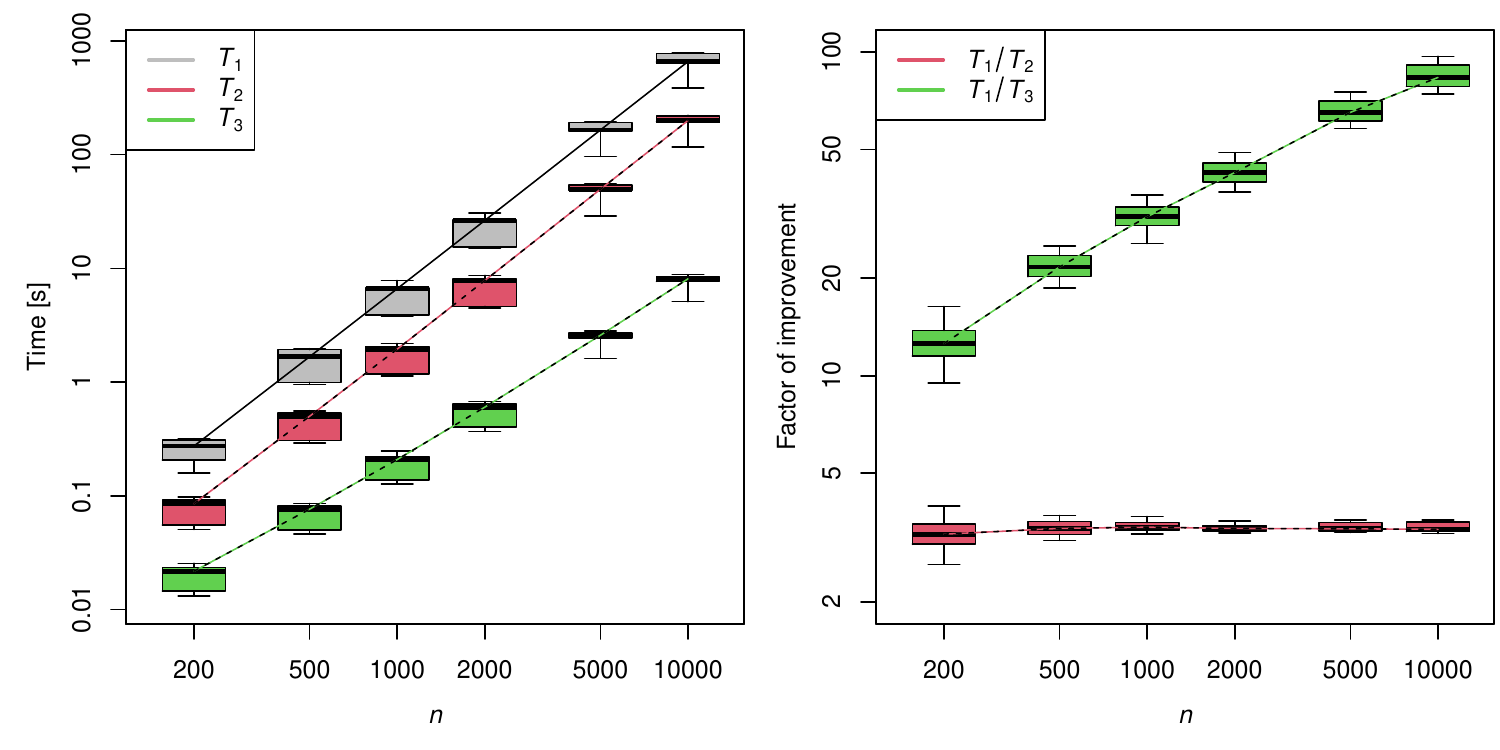}
\caption{Boxplots of computation times and ratios of running times for varying sample sizes. The whiskers indicate the $10\%$ and $90\%$ sample quantiles. The other elements of the boxplots are standard. A logarithmic scale was used for both axes.}
\label{fig:Boxplots}
\end{figure}

Figure~\ref{fig:Boxplots} displays the result of simulation experiments for sample sizes ranging from $200$ to $10\,000$, where the data were generated using the procedure mentioned earlier. The simulations indicate that the improvement due to using modified instead of standard PAVA is almost constant in $n$, whereas the improvement due to abridged instead of modified PAVA increases with $n$. Presumably, the complexity of the abridged PAVA for computing the isotonic distributional regression remains quadratic in $n$. But our numerical experiments show that the constant is substantially smaller than the one resulting from applying the usual PAVA with complexity $O(n)$ for $n - 1$ different levels of the response.

\paragraph{Numerical experiment 2.}
The goal of this experiment is to study the influence of the strength of the monotone association between $X$ and $Y$ on the efficiency gain of the abridged PAVA for isotonic distributional regression. The gains of abridged PAVA are expected to be milder when $Y$ is independent of $X$, and to become larger as the monotone association strengthens. The reason behind it is that, while the standard PAVA proceeds independently of the stochastic order, the abridged PAVA relies on the index $j_o$ indicating the component increasing in $\bs{z}(t-1)$ and on the nature of the partition corresponding to $A(\bs{z}(t-1))$, at a certain state $t \in \{1,\ldots,n\}$ of the procedure. If the monotone association is weak, then the partition corresponding to $A(\bs{z}(t-1))$ tends to contain fewer blocks in total and relatively large blocks in the middle of $\{1,\ldots,n\}$. If the index $j_o$ happens to lie in a block containing many indices to the right of $j_o$, even the abridged PAVA will have to inspect all of these.

To demonstrate this claim, we simulated $n$ independent bivariate Gaussian random vectors $(X,Y)^\top$ with correlation $\Corr(X,Y) = \rho \ge 0$. Note that the respective means and variances of $X$ and $Y$ have no influence on the results of the experiment. Indeed, the running times are invariant under strictly isotonic transformations of $X$ and of $Y$. In particular, the simulations for $\rho = 0$ cover all situations in which $X$ and $Y$ are stochastically independent with continuous distribution functions. The stochastic order between $\mathcal{L}(Y \vert X = x_1)$ and $\mathcal{L}(Y\vert X = x_2)$ for $x_1 < x_2$ becomes stronger as the correlation $\rho \in [0,1)$ increases, from an equality in distribution when $\rho = 0$ to a deterministic ordering when $\rho$ approaches $1$. Now, for sample sizes $n$ ranging from $200$ to $10\,000$ and for each correlation $\rho \in \{0, 0.5, 0.9\}$, the mean and standard deviation of the time ratio $T_3/T_1$ were estimated from $1\,000$ repetitions. The results are summarized in Table \ref{tab:GaussianSimulations}. As expected, the efficiency gain is smallest for $\rho = 0$. But even then, it is larger than $6$ for $n \ge 200$ and larger than $10$ for $n \ge 1\,000$.

\begin{table}
\[	
	\begin{array}{|r||rr||rr||rr|}
	\hline
	n &
	\multicolumn{2}{|c||}{\rho = 0}   &
	\multicolumn{2}{ c||}{\rho = 0.5} &
	\multicolumn{2}{ c| }{\rho = 0.9} \\
	\hline\hline
	200 &  6.5337 & (2.3496) & 10.7581 &  (3.6704) &  13.5695 &  (4.6390) \\
	\hline
    500 &  8.3029 & (2.6393) & 18.7010 &  (5.7806) &  26.1813 &  (8.0763) \\
	\hline
 1\,000 &  9.1351 & (3.1800) & 27.6290 &  (7.5007) &  41.4116 & (11.2161) \\
	\hline
 2\,000 &  9.7559 & (3.3532) & 39.3180 & (10.0337) &  62.8382 & (16.3293) \\
	\hline
 5\,000 & 10.7495 & (4.0525) & 62.4600 & (18.2002) & 108.4198 & (31.1414) \\
	\hline
10\,000 & 12.5190 & (5.6193) & 91.9084 & (33.4657) & 168.5030 & (58.7712) \\
	\hline
	\end{array}
\]
\caption{Means (and standard deviations) of the factor of improvement $T_3/T_1$ for different correlation values $\rho$ between $X$ and $Y$ and sample sizes $n$.}
\label{tab:GaussianSimulations}
\end{table}

\paragraph{Acknowledgments.}
The authors are grateful to a reviewer for constructive comments. This work was supported by Swiss National Science Foundation. R code is available at \url{https://github.com/AlexanderHenzi/abridgedPava}.


\end{document}